   \documentstyle[11pt,leqno]{article}
   \font\maxid=eufm10
   
   \def\p{\hbox{\maxid p}}
   \def\q{\hbox{\maxid q}}
   \def\Xg{\hbox{\maxid X}}

   \newtheorem{thm}{Theorem}[section]
   \newtheorem{lem}[thm]{Lemma}
   \newtheorem{pro}[thm]{Proposition}
   \newtheorem{cor}[thm]{Corollary}

   \newenvironment{dfn}
    {\refstepcounter{thm}
    {\medskip\par\noindent\bf Definition
    \arabic{section}.\arabic{thm} }}{\vspace{1ex}}
   \newenvironment{rem}
    {\refstepcounter{thm}
    {\medskip\par\noindent\bf Remark
    \arabic{section}.\arabic{thm} }}{\vspace{1ex}}
   \newenvironment{exa}
    {\refstepcounter{thm}
    {\medskip\par\noindent\bf Example
    \arabic{section}.\arabic{thm} }}{\vspace{1ex}}
   
   \def\longhookrightarrow{\lhook\joinrel\longrightarrow}

   \def\Bigdownarrow{\vphantom{\bigg|}\Big\downarrow}

   \def\cart{\lower2pt\hbox{\vbox{\hrule width9pt height0.4pt
             \hbox to9pt{\vrule height9pt width 0.4pt
             \hfil\vrule height9pt width 0.4pt}
             \hrule width9pt height 0.4pt}}}
   \def\Spec{\mathop{\mbox{\rm Spec}}\nolimits}
   \def\Ker{\mathop{\mbox{\rm Ker}}\nolimits}
   \def\Coker{\mathop{\mbox{\rm Coker}}\nolimits}
   
   \def\Pic{\mathop{\mbox{\rm Pic}}\nolimits}
   \def\Hom{\mathop{\mbox{\rm Hom}}\nolimits}
   
   \def\Extr{\mathop{\mbox{\rm Ext}}\nolimits}

   \def\gp#1{#1^{\mbox{\rm gp}}}

   \def\diag{\hbox{\maxid{d}}}
   \def\underrel#1#2{\mathrel{\mathop{#1}\limits_{#2}}}
   \newcommand{\pf}{{\sc Proof.}\hspace{2mm}}
   \newcommand{\qed}{$\Box$}
   \newcommand{\Na}{\mbox{\bf N}}
   \newcommand{\Z}{\mbox{\bf Z}}

   \renewcommand{\H}{\mbox{\rm H}}
   \renewcommand{\L}{{\cal L}}
   \newcommand{\M}{{\cal M}}
   
   \newcommand{\N}{{\cal N}}
   
   \renewcommand{\O}{{\cal O}}
   \newcommand{\Oi}{{\cal O}^\times}
   \newcommand{\T}{{\cal T}}

   \newcommand{\K}{{\cal K}}

   \newcommand{\la}{\lambda}
   \newcommand{\La}{\Lambda}
   \newcommand{\Ext}{{\cal E}xt}

   \title{\bf Logarithmic embeddings and logarithmic semistable reductions}
   \author{Fumiharu Kato\thanks{
   \hspace*{1.5em}{\em $1991$ Mathematics Subject Classification\/}.
   Primary 14B20;
   Secondary 13D10, 14D15, 16S80.
   }}
   \date{}
   
\begin{document}
   \def\footnotemark{\relax}
   \maketitle
   \begin{abstract}
   In this paper, we give a criterion for the existence of {\it logarithmic
   embeddings} --- which was first introduced by Steenbrink --- for
   general normal crossing varieties. Using this criterion, we also give a
   new proof of the theorem of Kawamata--Namikawa which gives
   a criterion for the existence of the log structure of semistable type.
   \end{abstract}
\section{Introduction}
Let $X$ be a connected, geometrically reduced algebraic scheme over a
field $k$. Then $X$ is said to be a {\it normal crossing variety} of
dimension $n-1$ if there exists an isomorphism of $k$-algebras
$$
\widehat{\O}_{X,x}\stackrel{\sim}{\longrightarrow}
k(x)[[T_1,\ldots,T_n]]/(T_1\cdots T_{l_x})
$$
for each closed point $x\in X$, where $\widehat{\O}_{X,x}$ denotes the
completion of the local ring $\O_{X,x}$ along its maximal ideal
(Definition \ref{ncvdef}).
Normal crossing varieties usually appear in contexts of algebraic geometry
via degenerations and normal crossing divisors.
In the first case, they appear as a specialization of a family of smooth
varieties.
Normal crossing varieties are usually considered and expected to be limits
of smooth varieties, and --- as is well--known --- they
are important to the theory of moduli.
As for the second situation, a {\it normal crossing divisor} is a divisor
of a smooth variety which itself is a normal crossing variety.
Normal crossing divisors play important roles
in various fields of algebraic geometry.
For example, a pair of smooth variety and its normal crossing divisor is
usually called a log variety.
Considering log varieties instead of smooth varieties ---
or usually admitting some mild singularities --- alone,
some algebro geometric theories ({\it e.g.}, minimal model theory, etc.)
are well generalized.

Relating with a normal crossing variety $X$, there are two problems,
{\it smoothings} and {\it embeddings},
in light of degenerations and normal crossing divisors, respectively.

The {\it smoothing problem} is a problem to find a Cartesian diagram
$$
\begin{array}{ccc}
X&\longrightarrow&\Xg\\
\Bigdownarrow&&\Bigdownarrow\\
0&\longrightarrow&\Delta\rlap{,}
\end{array}
$$
for a normal crossing variety $X$ (in this situation, we should assume that
$X$ is proper over $k$), where $\Delta$ is a one-dimensional regular scheme,
$\Xg$ is a regular scheme proper flat and generically smooth over $\Delta$,
and $0$ is a closed point of $\Delta$ whose residue field is $k$.
We usually take, as the base scheme $\Delta$, the spectrum of a discrete
valuation ring, {\it e.g.}, the
ring of formal power series over $k$ or --- in case $k$ is perfect ---
the ring of Witt vectors over $k$.
In the complex analytic situation, Friedmann \cite{Fri1}
studied the smoothing problem generally, and solve it
for degenerated K3 surfaces.
Recently, Kawamata--Namikawa \cite{K-N1} approached this problem
by introducing a new method; the {\it logarithmic}
method.
The Cartesian diagram as above with $\Delta$ a
spectrum of an Artinian local ring $A$ is
called an {\it infinitesimal smoothing},
if it is \'{e}tale locally isomorphic to the diagram
$$
\begin{array}{ccc}
\Spec k[Z_1,\ldots,Z_n]/(Z_1\cdots Z_l)&\longrightarrow&
\Spec A[Z_1,\ldots,Z_n]/(Z_1\cdots Z_l-\pi)\\
\Bigdownarrow&&\Bigdownarrow\\
\Spec k&\longrightarrow&\Spec A\rlap{,}
\end{array}
$$
where $\pi$ is an element of the maximal ideal of $A$ and $\pi\neq 0$.
The central problem to find such an infinitesimal smoothing is to compute
the obstruction class of $X$ to have such a diagram and to show vanishing
or non--vanishing of it.

The {\it embedding problem} is a problem to find
a closed embedding $X\hookrightarrow V$ over $k$ of $X$ as a normal crossing
divisor, where $V$ is a smooth variety over $k$.
If $X$ is smoothable in the above sense
with $\Delta$ a smooth algebraic variety over $\Spec k$, the smoothing family
$X\hookrightarrow\Xg$ gives an embedding of this sense.
If $X$ is smooth, this problem becomes trivial, since we can take as $V$ the
product of $X$ and, for example, ${\mbox{\bf P}}^1$.
But for a general normal crossing variety, this problem seems
far from satisfactory solutions.
Similarly to the smoothing problem, we can consider
this problem in the infinitesimal sense.

In this paper, we consider the above problems in a {\it logarithmic} sense.
We consider logarithmic generalizations of smoothings and embeddings of
normal crossing varieties according to Kajiwara \cite{Kaj1},
Kawamata--Namikawa \cite{K-N1} and Steenbrink \cite{Ste1}, and we solve
their existence problems.
These generalizations are done in terms of logarithmic geometry of
Fontaine, Illusie and Kazuya Kato.

{\it Logarithmic geometry} --- or {\it log geometry} --- was first founded by
Fontaine and Illusie  based on their idea of, so--called,
{\it log structures};
afterwards, it was established as a generally
organized theory and applied to various fields of algebraic and
arithmetic geometry by Fontaine, Illusie and Kazuya Kato (cf. \cite{Kat1},
\cite{Kat4}).
In various kinds of geometries including algebraic geometry,
we usually consider local ringed spaces, {\it i.e.}, the pairs of topological
spaces --- possibly in the sense of Grothendieck topologies ---
and sheaves of local rings over them.
The basic idea of Fontaine and Illusie is that, instead of
local ringed spaces alone, they consider local ringed spaces equipped with some
additional structure --- which they call the logarithmic structures ---
written in
terms of sheaves of commutative and unitary monoids (see \cite{Kat1} for the
precise definition). In algebro geometric situations,
these log structures usually
represent ``something'' of the underlying local ringed spaces, {\it e.g.},
divisors or the structure of torus embeddings, etc. Through these
foundations, they suggested to generalize the ``classical'' geometries by
considering ``log objects'' --- such as {\it log schemes} --- which are
the pairs of local ringed spaces and log structures on them.

In the present paper, we recall and generalize the
{\it logarithmic embedding}
(Definition \ref{logembdef}) introduced by Steenbrink \cite{Ste1}.
A logarithmic embedding --- which is regarded as a logarithmic generalization
of a log variety --- is a certain log scheme $(X,\M_X)$ with $X$ a
normal crossing variety.
Then we prove the following theorem which
gives a criterion for the existence of logarithmic embeddings:

\vspace{3mm}\noindent
{\bf Theorem}\ ({\it Theorem \ref{mainthm}})\ {\it
For a normal crossing variety $X$, a logarithmic embedding of $X$ exists if
and only if there exists a line bundle ${\cal L}$ on $X$ such that
${\cal L}\otimes_{\O_X}\O_D\stackrel{\sim}{\rightarrow}
\T^1_X$, where $D$ is the singular locus of $X$.}

\vspace{3mm}
Here, $\T^1_X$ is an invertible $\O_D$-module, called the {\it infinitesimal
normal bundle} (cf. \cite{Fri1}), which is naturally isomorphic to
$\Ext^1_{\O_X}(\Omega^1_X,\O_X)$; we recall the construction of it in
\S 3.

A normal crossing variety $X$ is said to be
{\it $d$-semistable} if
$\T^1_X$ is a trivial bundle on $D$ (cf. \cite{Fri1}).
By the above theorem, any $d$-semistable normal crossing variety $X$ has a
logarithmic embedding.

As for the smoothing problem, we recall and generalize the concept,
the {\it logarithmic semistable reduction} (Definition \ref{logsemidef})
introduced by Kajiwara \cite{Kaj1} (in one dimensional case) and
Kawamata--Namikawa \cite{K-N1} (by a different but essentially the same
method).
Using the above theorem, we get a criterion for the existence of
logarithmic semistable reductions, which was first proved by
Kawamata--Namikawa \cite{K-N1} in the complex analytic situation, as
follows:

\vspace{3mm}\noindent
{\bf Theorem}\ ({\it Theorem \ref{mainthm2}})
{\rm (cf. \cite{K-N1})}\
{\it For a normal crossing variety $X$,
the log structure of semistable type
on $X$ exists if and only if $X$ is $d$-semistable.}

\vspace{3mm}
The composition of this paper is as follows.
In \S 2, we study the geometry of normal crossing varieties in general.
In particular, we define good \'{e}tale local charts on normal crossing
varieties, and prove the existence of them.
In \S 3, we recall the basic construction of the tangent complex of a
normal crossing variety, and introduce the invertible sheaf
${\cal T}^1_X$ on $D$.
We introduce the logarithmic embedding in \S 4.
This section also contains the proof of our main theorem.
The logarithmic semistable reduction is studied in \S 5.

The author thanks T. Fujisawa for useful communications.
The author is also grateful to Professors K. Ueno, S. Usui and T. Yusa
for their helpful comments.

{\sc Conventions}:\
All sheaves are considered with respect to \'{e}tale topology.
By a monoid, we mean --- as usual in the contexts of log geometry ---
a set with a commutative and associative binary operation and the neutral
element.
For such a monoid $M$, we denote by $\gp{M}$ the Grothendieck group of $M$.
We denote by $\Na$ the monoid of non--negative integers.

\section{Normal crossing varieties}
Throughout this paper, we always work over a fixed base field $k$.
As usual, an algebraic $k$-scheme is, by definition, a seperated
scheme of finite type over $k$.
Let $X$ be an algebraic $k$-scheme and $x\in X$ a point.
We denote the residue field at $x\in X$ by $k(x)$.
\begin{dfn}\label{ncvdef}
Let $X$ be a connected and geometrically reduced algebraic $k$-scheme.
Then $X$ is said to be a {\it normal crossing variety} over $k$ of
dimension $n-1$ if the following condition is satisfied: For any closed
point $x\in X$, there exists an isomorphism
\begin{equation}\label{ncvdefloc}
\widehat{\O}_{X,x}\stackrel{\sim}{\longrightarrow}
k(x)[[T_1,\ldots,T_n]]/(T_1\cdots T_{l_x})
\end{equation}
of $k$-algebras, where $l_x$ is an integer $(1\leq l_x\leq n)$
depending on $x$.
Here, we denote by $\widehat{\O}_{X,x}$ the completion of the local ring
$\O_{X,x}$ by its maximal ideal.
\end{dfn}

The integer $l_x$ is called the {\it multiplicity} at $x\in X$.
We sometimes denote it by $l^X_{x}$ if we want to emphasize the scheme $X$.
The Zariski closure of the set of closed points whose multiplicity is greater
than 1 is the singular locus of $X$, which we denote by $D$.

A standard example of normal crossing varieties is an affine scheme
\begin{equation}\label{ncvstandard}
\Spec k[T_1,\ldots,T_n]/(T_1\cdots T_l)\ \ (1\leq l\leq n).
\end{equation}
This scheme consists of $l$ irreducible components which intersect
transversally along the singular locus
\begin{equation}
\Spec k[T_1,\ldots,T_n]/(T_1\cdots\widehat{T_j}\cdots T_l\ :\ 1\leq j\leq l).
\end{equation}
Each irreducible component is isomorphic to the affine $(n-1)$-space
over $k$.
In general, a normal crossing variety $X$ is said to be {\it simple}
if each irreducible component of $X$ is smooth over $k$.
For example, a smooth $k$-variety is a simple normal crossing variety.

Let $V$ be a smooth $k$-variety of dimension $n$.
A reduced divisor $X$ on $V$ is called a {\it normal crossing divisor}
if $X$ itself is a normal crossing variety
of dimension $n-1$.
In this case, the closed embedding $X\hookrightarrow V$ is called a {\it NCD
embedding} of $X$.
For example, the affine normal crossing variety (\ref{ncvstandard}) is
a normal crossing divisor in the affine $n$-space over $k$.

The proof of the following proposition is straightforward and is left
to the reader.
\begin{pro}\label{ncvetale}
Let $Y$ be a connected scheme \'{e}tale over a connected algebraic
$k$-scheme $X$.
If $X$ is a normal crossing variety, then so is $Y$.
The converse is also true if the \'{e}tale morphism $Y\rightarrow X$
is surjective.
\end{pro}

It is clear that an \'{e}tale morphism leaves invariant the multiplicity at
every closed point, {\it i.e.}, if $\varphi\colon Y\rightarrow X$ is an
\'{e}tale
morphism of normal crossing varieties and $y\in Y$ is a closed point,
then we have $l^Y_y=l^X_{\varphi(y)}$.

In the following paragraphs of this section, we shall study the local nature
of normal crossing varieties for the later purpose.
In the subsequent sections, we need to take a good \'{e}tale
neighborhood around every closed point.
We require that these \'{e}tale neighborhoods have good coordinate systems
which serve for several explicit calculations.
To clarify the notion of ``good'' \'{e}tale neighborhoods, we define them
as follows:
\begin{dfn}\label{ncvchart}{\rm
Let $X$ be a normal crossing variety and $x\in X$ a closed point.
Let $\varphi\colon U\rightarrow X$ be an \'{e}tale morphism with $U$ a simple
normal crossing variety and
$z_1,\ldots,z_{l_x}\in\Gamma(U,\O_U)$, where $l_x$ is the multiplicity
at $x$.
Then $(\varphi\colon U\rightarrow X;\ z_1,\ldots,z_{l_x})$ is said to be a {\it
local chart} around $x$ if the following conditions are satisfied:
\begin{description}
\item[{\rm (a)}] There exists a unique point $y\in U$ such that
$\varphi(y)=x$.
\item[{\rm (b)}] There exists a closed immersion $\iota\colon U\hookrightarrow
V$, where $V$ is an affine smooth $k$-scheme.
\item[{\rm (c)}] There exist $Z_1,\ldots,Z_n\in\Gamma(V,\O_V)$ which form
a regular parameter system at $\iota(y)\in V$ such that $z_i=\iota^*Z_i$
for $1\leq i\leq l_x$, and $U$ is defined as a closed
subset in $V$ by the ideal $(Z_1\cdots Z_{l_x})$.
\item[{\rm (d)}] each ideal
$(z_i)$ is prime and the irreducible components of $U$ are precisely the
closed subsets of $U$ corresponding to the ideals $(z_1),\ldots,(z_{l_x})$.
\end{description}
}
\end{dfn}

Note that $\iota\colon U\hookrightarrow V$ is, due to (c), a NCD embedding.
Moreover, due to (d), all the irreducible components intersect and contain
the point $y$.

The following theorem assures the existence of local chart around every
closed point of normal crossing variety $X$.
We prove this theorem later in this section.
\begin{thm}\label{chartexist}
Let $X$ be a normal crossing variety and $x\in X$ a closed point.
Then there exists a local chart
$(\varphi\colon U\rightarrow X;\ z_1,\ldots,z_{l_x})$ around $x$.
\end{thm}

Since any \'{e}tale open set of $X$ is again a normal crossing variety,
we have the following:
\begin{cor}\label{specialcov}
Let $X$ be a normal crossing variety.
Then the set of all local charts forms an open basis with respect to the
\'{e}tale topology on $X$.
\end{cor}

\begin{rem}\label{ncvemb}{\rm
Theorem \ref{chartexist} implies that any normal crossing variety is
realized as a simple normal crossing divisor on some smooth $k$-variety
\'{e}tale locally.
But a normal crossing variety, in general, cannot be a normal crossing
divisor globally on a smooth $k$-variety.
In the next section, we will see a necessary condition for a normal
crossing variety to be a normal crossing divisor
(Proposition \ref{suffcondemb}).
}
\end{rem}

For the proof of Theorem \ref{chartexist}, we need one lemma:
\begin{lem}\label{htzero}
Let $\q$ be a height zero prime ideal in $K[T_1,\ldots,T_n]/(T_1\cdots T_l)$
$(1\leq l\leq n)$, where $K$ is a field.
Then $\q=(T_j)$ for some $j$ $(1\leq j\leq l)$.
\end{lem}
\pf
By Krull's principal ideal theorem, any non--zero element in $\q$ is a
zero factor. Hence any element in $\q$ is a multiple of $T_j$'s $(1\leq
j\leq l)$. Since $\q$ is a prime ideal, $\q$ must contain $T_j$ for some
$j$ $(1\leq j\leq l)$, {\it i.e.}, $(T_j)\subseteq\q$. But since the height of
$\q$ is zero and $(T_j)$ is a prime ideal, we have $\q=(T_j)$.
\qed

\vspace{3mm}
{\sc Proof of Theorem \ref{chartexist}.}\hspace{2mm}
The complete local ring $\widehat{\O}_{X,x}$ is isomorphic to the complete
local
ring $k(x)[[T_1,\ldots,T_n]]/(T_1\cdots T_{l_x})$ which is a completion of
the local ring $(k(x)[T_1,\ldots,T_n]/(T_1\cdots T_{l_x}))_{0}$.
Then due to \cite[Corollary (2.6)]{Art1}, there exist a scheme
$U$ and \'{e}tale morphisms $\varphi\colon U\rightarrow X$ and
$\phi\colon U\rightarrow\Spec k(x)[T_1,\ldots,T_n]/(T_1\cdots T_{l_x})$ such
that
$\varphi(y)=x$ and $\phi(y)=0$ for some $y\in U$.
We fix this closed point $y\in U$.
Since $\varphi$ is \'{e}tale, we may assume --- replacing $U$ by its Zariski
open subset if necessary --- that $y$ is the only point which
is mapped to $x$ by $\varphi$.
Obviously we may assume that $U$ is connected and affine.
We can remove all the irreducible components which do not contain $y$.
Then we may assume that all the irreducible components of $U$ contain $y$.
We set $U=\Spec A$ and $B\colon =k(x)[T_1,\ldots,T_n]/(T_1\cdots T_{l_x})$.

Since $U$ is \'{e}tale over a reduced $k$-scheme $\Spec B$,
the $k$-algebra $A$ is reduced.
Take a minimal prime factorization
\begin{equation}\label{prmfac}
(0)=\p_1\cap\cdots\cap\p_{l}.
\end{equation}
of the ideal $(0)=\sqrt{(0)}$.
Since each $\p_i$ is minimal in the set of all prime ideals, the height of
each $\p_i$ is zero.
Obviously the prime decomposition (\ref{prmfac}) precisely corresponds to the
decomposition of $U$ into irreducible components.
Set $\q_i\colon =\phi(\p_i)$ which is a prime ideal of height zero in $B$ for
$1\leq i\leq l$.
Due to Lemma \ref{htzero}, we have $\q_i=(T_{j_i})$ for some $j_i$
$(1\leq j_i\leq l_x)$, {\it i.e.}, any generic point of a irreducible component
of
$U$ is mapped by $\phi$ to a generic point of a irreducible component of
$\Spec B$.

Let us suppose that the map $i\mapsto j_i$ is not injective,
{\it i.e.}, there exist $i$ and $j$
$(i\neq j)$ such that $\q_i=\q_j$.
Consider the Cartesian diagram
$$
\begin{array}{ccc}
\overline{\{\q_i\}}\times_{\Spec B}U&\longhookrightarrow&U\\
\llap{$\phi_i$}\Bigdownarrow&&\Bigdownarrow\rlap{$\phi$}\\
\overline{\{\q_i\}}&\longhookrightarrow&\Spec B\rlap{,}
\end{array}
$$
where the horizontal arrows are closed immersions and the vertical ones are
\'{e}tale.
The scheme $U_i\colon =\overline{\{\q_i\}}\times_{\Spec B}U$ is also a normal
crossing variety. Since $\overline{\{\p_i\}}\cap\overline{\{\p_j\}}$ is a
closed subscheme (which contains $y$) of $U_i$,
the multiplicity $l^{U_i}_y$ at $y$ in $U_i$ is greater than 1.
But since the irreducible component $\overline{\{\q_i\}}$ is smooth,
we have $l_{\phi_i(y)}=1$.
This is a contradiction since $l^{U_i}_y=l_{\phi_i(y)}$.
Thus, the map $i\mapsto j_i$ is injective, {\it i.e.}, there is at
most one component over each component of $\Spec B$.
Moreover, in this case, we have $U_i=\overline{\{\p_i\}}$.
Then the irreducible component $U_i$ is \'{e}tale
over a smooth scheme $\overline{\{\q_i\}}$,
and hence the normal crossing
variety $U$ is simple. Moreover, the prime ideal $\p_i$ is a principal ideal
$(z_i)$, where $z_i\colon =\phi^*T_{j_i}$ for $1\leq i\leq l$,
since $\overline{\{\p_i\}}=\overline{\{\q_i\}}\times_{\Spec B}U$ implies that
$\p_i=\q_i\otimes_BA$.

Since the map $i\mapsto j_i$ is injective, we have $l\leq l_x$.
Note that the multiplicity $l_y$ at y in $U$ equals to $l_x$.
Since the simple normal crossing variety $U$ consists of $l$ irreducible
components, we have $l_y=l_x\leq l$.
Hence we have $l_x=l$.

The scheme $\Spec B$ is a normal crossing divisor in the $n$ dimensional
affine space over $k(x)$.
Hence, due to \cite[Expos\'{e} 1. Proposition 8.1]{Gro1}, any point in $U$
has a Zariski open neighborhood which is embedded in a smooth
$k(x)$-variety as a normal crossing divisor.
This implies that, replacing $U$ by its Zariski open neighborhood of $y$,
we may assume that $U$ can be embedded in an affine smooth $k(x)$-scheme
$V=\Spec R$ of dimension $n$ as a normal crossing divisor.
Let $\iota\colon U\hookrightarrow V$ be the closed immersion.

Finally, consider the Cartesian diagram
$$
\begin{array}{ccc}
U&\stackrel{\iota}{\longhookrightarrow}&V\\
\llap{$\phi$}\Bigdownarrow&&\Bigdownarrow\rlap{$\Phi$}\\
\Spec B&\longhookrightarrow&\Spec k[T_1,\ldots,T_n]\rlap{,}
\end{array}
$$
where $\Phi$ is an \'{e}tale morphism.
Set $Z_i\colon =\Phi^*T_i\in\Gamma(V,\O_V)$ for $1\leq i\leq n$.
Then $Z_1,\ldots,Z_n$ form a regular parameter system at $\iota(y)\in V$.
We also have $z_i=\iota^* Z_i$ for $1\leq i\leq l_x$.
It is clear that the closed subscheme $U$ in $V$ is defined by an ideal
$(Z_1\cdots Z_{l_x})$.
Then the proof of the theorem is completed.
\qed

The following lemma will be needed in the later arguments.
\begin{lem}\label{inj-zero}
Let $(\varphi'\colon U'\rightarrow X;z'_1,\ldots,z'_{l'})$ be a local chart
on $X$ around some closed point and
$(\psi\colon U\rightarrow U';z_1,\ldots,z_{l})$ a local chart on $U'$
around some closed point.
Then $\psi\colon U\rightarrow U'$ is injective in codimension zero,
{\it i.e.}, it maps the generic points of irreducible components on $U$
injectively to those of $U'$.
\end{lem}

\pf
Let $\eta\in U'$ be a codimension zero point.
Since $U'$ is simple, $\overline{\{\eta\}}$ is regular and so is
$\overline{\{\eta\}}\times_{U'}U$ whenever it is not empty.
Then each connected component of $\overline{\{\eta\}}\times_{U'}U$ is
irreducible and its generic point is of codimension zero.
Hence each connected component of $\overline{\{\eta\}}\times_{U'}U$ is
an irreducible component of $U$. Since any two of irreducible
components of $U$ intersect, $\overline{\{\eta\}}\times_{U'}U$ itself is an
irreducible component of $U$.
Hence, if $\xi\in U$ is a codimension zero point such that $\psi(\xi)=\eta$,
we have $\overline{\{\xi\}}=\overline{\{\eta\}}\times_{U'}U$.
In particular, there exists at most one such $\xi$.
\qed

For a normal crossing variety $X$, the {\it normalization}
$\nu\colon\widetilde{X}\rightarrow X$ of $X$ is defined as usual:
The scheme $\widetilde{X}$ is defined by the disjont union of the
normalizations of irreducible components of $X$ and $\nu\colon\widetilde{X}
\rightarrow X$ is the natural morphism. The normalization
$\widetilde{X}$ is a smooth $k$-scheme due to Theorem \ref{chartexist} and
the following lemma.

\begin{lem}\label{normalization}
Let $U\rightarrow Z$ be a \'{e}tale morphism of $k$-varieties.
Let $\widetilde{U}\rightarrow U$ and $\widetilde{Z}\rightarrow Z$ be
normalizations of $U$ and $Z$, respectively. Then there exists a natural
isomorphism
$\widetilde{U}\stackrel{\sim}{\rightarrow}U\times_{Z}\widetilde{Z}$.
In particular, the natural morphism $\widetilde{U}\rightarrow\widetilde{Z}$
is \'{e}tale.
\end{lem}

\pf
Since $U\times_{Z}\widetilde{Z}\rightarrow\widetilde{Z}$ is \'{e}tale and
$\widetilde{Z}$ is normal, $U\times_{Z}\widetilde{Z}$ is a normal variety.
Hence there exists a unique morphism $\phi\colon
U\times_{Z}\widetilde{Z}\rightarrow
\widetilde{U}$ which factors the morphism $U\times_{Z}\widetilde{Z}\rightarrow
U$. Moreover $\phi$ also factors the morphism $U\times_{Z}\widetilde{Z}
\rightarrow\widetilde{Z}$ since the last morphism is the unique morphism
determined by the morphism $U\times_{Z}\widetilde{Z}\rightarrow Z$. Hence
the natural morphism $\varphi\colon \widetilde{U}\rightarrow
U\times_{Z}\widetilde{Z}$
is the inverse morphism of $\phi$.
\qed

For a local chart $(\varphi\colon U\rightarrow X; z_1,\ldots,z_l)$, the
normalization of $U$ is given by the disjoint union of all irreducible
components and the natural morphism, {\it i.e.},
$$
\nu_U\colon\widetilde{U}=\coprod^l_{i=1}U_i\longrightarrow U,
$$
where $U_i$ is the irreducible component of $U$ corresponding to the
ideal $(z_i)$.

Set $\overline{D}\colon =D\times_X\widetilde{X}$, which is a divisor of
$\widetilde{X}$.
\begin{lem}
$\overline{D}$ is a normal crossing divisor of $\widetilde{X}$.
\end{lem}

\pf
Let $(\varphi\colon U\rightarrow X; z_1,\ldots,z_l)$ be a local chart on $X$.
Then $D_U\colon =D\times_{X}U$ is nothing but the singular locus of $U$ and is
\'{e}tale over $D$.
Consider the normalization $\nu_U\colon \widetilde{U}\rightarrow U$ as above.
Set
$$
\overline{D_U}\colon =D_U\times_{U}\widetilde{U}.
$$
Clearly, $\overline{D_U}$ is a normal crossing divisor of $\widetilde{U}$
defined by an ideal $(z_1\cdots\widehat{z_i}\cdots z_l)$ on $U_i$.
There exists a natural morphism $\overline{D_U}\rightarrow\overline{D}$.
Since one can easily see that there exists a natural isomorphism
$$
\overline{D_U}\stackrel{\sim}{\longrightarrow}
\overline{D}\times_{\widetilde{X}}\widetilde{U}
$$
and the morphism $\widetilde{U}\rightarrow\widetilde{X}$ is \'{e}tale
due to Lemma \ref{normalization}, the morphism $\overline{D_U}\rightarrow
\overline{D}$ is \'{e}tale.
Then, considering all the local charts on $X$, $\overline{D}$ is a normal
crossing divisor on $\widetilde{X}$ due to Proposition \ref{ncvetale}.
\qed

\section{Tangent complex on a normal crossing variety}
In this section, we recall the tangent complex and the infinitesimal normal
bundle $\T^1_X$ of a normal crossing variety $X$ which will play
important roles in the subsequent sections.

Let $X$ be a normal crossing variety over a field $k$.
For a local chart $(\varphi\colon U=\Spec A\rightarrow X; z_1,\ldots,z_l)$ of
$X$ around some closed point, we use the folowing notation in this and
subsequent sections:
Let $V=\Spec R$ and $Z_1,\ldots,Z_l$ be as in Definition \ref{ncvchart}.
Set $I_j\colon=(Z_j)$ and $J_j\colon=(Z_1\cdots\widehat{Z_j}\cdots Z_l)$
for $1\leq j\leq l$. (If $l=1$, we set $J_1=R$ for the convention.)
Then $A=R/I$ where $I\colon=I_1\cdots I_l$.
Moreover, the ideal $I_j/I\subset A$ is generated by $z_j=(Z_j\,\mbox{\rm
mod}\,I)$
and is prime of height zero.
Set $J\colon=J_1+\cdots+J_l$.
Then the singular locus $D_U\colon=D\times_XU$ of $U$ is the closed subscheme
defined by $J$. We set $Q\colon=R/J$.
Note that, for $1\leq j\leq l$, $I_j/II_j$ is a free $A$-module of rank
one and is generated by $\zeta_j\colon=(Z_j\,\mbox{\rm mod}\,II_j)$.
There exists a natural isomorphism $I_j/II_j\otimes_AQ\stackrel{\sim}
{\rightarrow}I_j/JI_j$ of $Q$-modules which maps $\zeta_j\otimes 1$ to
$\xi_j\colon=(Z_j\,\mbox{\rm mod}\,JI_j)$.
Moreover, there exists a natural isomorphism
\begin{equation}\label{conormal}
I/I^2\stackrel{\sim}{\rightarrow}
I_1/II_1\otimes_A\cdots\otimes_AI_l/II_l
\end{equation}
of $A$-modules, and hence, the $A$-module $I/I^2$ is free of rank one and is
generated by $\zeta_1\otimes\cdots\otimes\zeta_l$.
We denote by $\pi_j$ the natural projection $I_j/II_j\rightarrow I_j/I\subset
A$.

The cotangent complex of the morphism $k\rightarrow A$ is given by
$$
L^{\cdot}:0\longrightarrow R\otimes_{R}A\stackrel{\delta}{\longrightarrow}
\Omega^1_{R/k}\otimes_{R}A\longrightarrow 0,
$$
where $\delta$ is defined by $R\rightarrow F\cdot R
\stackrel{d}{\rightarrow}\Omega^1_{R/k}$ with $F\colon=Z_1\cdots Z_l$
(cf. \cite{L-S1}).
Then the tangent complex of $U$ is the complex
$$
\Hom_A(L^{\cdot},A):0\longrightarrow\Theta_{R/k}\otimes_{R}A
\stackrel{\delta^*}{\longrightarrow}\Hom_A(R\otimes_{R}A,A)
\longrightarrow 0,
$$
where $\Theta_{R/k}\colon =\Hom_R(\Omega^1_{R/k},R)$.

We define
\begin{equation}\label{tangentdef}
T^1_A=\Hom_A(R\otimes_{R}A,A)/\delta^*(\Theta_{R/k}\otimes_{R}A).
\end{equation}

\begin{lem}\label{tangentloc}
We have the natural isomorphism
\begin{equation}\label{tangentloc1}
T^1_A\stackrel{\sim}{\rightarrow}\Hom_A(I/I^2,A)\otimes_{A}Q.
\end{equation}
\end{lem}

\pf
Consider the exact sequence
$$
0\longrightarrow I/I^2\longrightarrow\Omega^1_{R/k}\otimes_{R}A
\longrightarrow\Omega^1_{A/k}\longrightarrow 0.
$$
By definition, we have
$T^1_A=\Coker(\Hom_A(\Omega^1_{R/k}\otimes_{R}A,A)\rightarrow
\Hom_A(I/I^2,A))$.
Then one can show --- by direct calculations --- that
$\Hom_A(I/I^2,A)\rightarrow T^1_A$ is nothing but the ``tensoring'' morphism
$\otimes_AQ$. Moreover, we have
$T^1_A\stackrel{\sim}{\rightarrow}\Extr^1_A(\Omega^1_{A/k},A)$.
\qed

Considering all the local charts $U$ on $X$, these modules $T^1_A$ glue to an
invertible $\O_D$-module on $X$, which is denoted by $\T^1_X$;
this is well--known (cf. \cite{L-S1}) but, for the later purpose, we prove
it in the following.

Suppose we have two local charts $(\varphi\colon U\rightarrow X;
z_1,\ldots,z_{l})$ and  $(\varphi'\colon U'\rightarrow X;
z'_1,\ldots,z'_{l'})$ and an \'{e}tale morphism $\psi\colon U\rightarrow U'$
such that $\varphi=\varphi'\circ\psi$.
(Because we are interested in the singular locus, we shall assume
$l>1$ and $l'>1$.)
For these local charts, we use all the notation as above. (For $U'$, we denote
them by $A'$, $I'$, $J'$, $\zeta'_j$, etc.)
Let $f\colon A'\rightarrow A$ be the ring homomorphism corresponding to
$\psi$.
We shall show that the morphism $\psi$ induces naturally an isomorphism
$T^1_{A'}\otimes_{Q'}Q\stackrel{\sim}{\rightarrow}T^1_A$ of $Q$-modules.

Let $U_j$ (resp. $U'_j$) be the irreducible component of $U$ (resp. $U'$)
corresponding to $I_j/I$ (resp. $I'_j/I'$) for $1\leq j\leq l$
(resp. $1\leq j\leq l'$).
Since $\psi$ is \'{e}tale and injective in codimension zero (Lemma
\ref{inj-zero}), we may assume that the generic point of $U_j$ is mapped to
that of $U'_j$ by $\psi$ for $1\leq j\leq l$.
In particular, we have $l\leq l'$.
Then one sees easily that $U\times_{U'}U'_j\cong U_j$ for $1\leq j\leq l$.
This implies that $A/(I_j/I)\cong (A'/(I'_j/I')\otimes_{A'}A)
(\cong A/((I'_j/I')\otimes_{A'}A))$, and hence,
\begin{equation}\label{ideals}
I_j/I=(I'_j/I')\otimes_{A'}A,\ (1\leq j\leq l)
\end{equation}
as ideals in $A$.
For $1\leq j\leq l$, we can set $f(z'_j)=u_jz_j$ for some $u_j\in A$.
Here, each $u_j$ is determined up to modulo $J_j/I$.
Due to (\ref{ideals}), $u_jz_j$ generates the ideal $I_j/I$,
and hence, $u_j$ is a unit in $A/(J_j/I)$ (and, of course, in $A/(J/I)$).
(Note that $u_j$ is not necessarily a unit in $A$, since $A$ is not an integral
domain for $l>1$.)
Then there exists an isomorphism (naturally induced by $f$) of $Q$-modules
\begin{equation}\label{def-tau}
\tau_j\colon I'_j/I'I'_j\otimes_{A'}Q\stackrel{\sim}{\longrightarrow}
I_j/II_j\otimes_AQ
\end{equation}
by $\xi'_j\mapsto(u_j\,\mbox{\rm mod}\,J/I)\xi_j$.
The natural projection $\pi'_j\colon I'_i/I'I'_i\rightarrow I'_i/I'
\subset A'$
$(1\leq i\leq l')$ and $f$ induce an $A$-module morphism
\begin{equation}\label{def-proj}
\widetilde{\rho}_i\colon I'_i/I'I'_i\otimes_{A'}A\longrightarrow A.
\end{equation}
For $1\leq j\leq l$, $\widetilde{\rho}_j$ maps $I'_i/I'I'_i\otimes_{A'}A$
surjectively onto
$I_j/I$, and for $i>l$, $\widetilde{\rho}_i$ is an isomorphism;
because, for $i>l$, one sees that
$\widetilde{\rho}_i(\zeta'_i\otimes 1)=f(z'_i)$
is an invertible element of $A$ as folows: Since $\psi$ is injective in
codimension zero, the point $I'_i/I'$ does not belong to $\psi(U)$;
hence $\psi$ maps $U=\Spec A$ to $\Spec A'_{(I'_i/I')}$, and this
implies the image of elements in $I'_i/I'$ under $f$ is invertible.
Set $\rho_i\colon=\widetilde{\rho}_i\otimes_AQ$.
Then these isomorphisms induce
\begin{equation}\label{def-tau2}
\tau\colon=\tau_1\otimes_Q\cdots\otimes_Q\tau_l\otimes_Q
\rho_{l+1}\otimes_Q\cdots\otimes_Q\rho_{l'}\colon
I'/I'^2\otimes_{A'}Q\stackrel{\sim}{\rightarrow}I/I^2\otimes_AQ.
\end{equation}
The $Q$-dual of $\tau$ is the desired isomorphism (cf. Lemma \ref{tangentloc}).
One can easily check that this isomorphism $\tau$ does not depends on
parameters $z'_j$, $z_j$; it is cannonically induced by
$f\colon A'\rightarrow A$.
Hence, for any sequence of \'{e}tale morphisms of local charts
$U\stackrel{\psi}{\rightarrow}U'\stackrel{\psi'}{\rightarrow}U''$, we obviously
have $\tau''=\tau\circ(\tau'\otimes_{Q'}Q)$, where
$\tau\colon I'/I'^2\otimes_{A'}Q\stackrel{\sim}{\rightarrow}I/I^2\otimes_AQ$,
$\tau'\colon I''/I''^2\otimes_{A''}Q'\stackrel{\sim}{\rightarrow}
I'/I'^2\otimes_{A'}Q'$ and
$\tau''\colon I''/I''^2\otimes_{A''}Q\stackrel{\sim}{\rightarrow}
I/I^2\otimes_AQ$ are the isomorphisms defined as above with respect to
$\psi$, $\psi'$ and $\psi'\circ\psi$, respectively.
Then one sees easily that there exists a unique $\O_D$-module whose restriction
to each $U$ is the $\O_{D_U}$-module corresponding to $T^1_A$; and it is
nothing but our desired $\O_D$-module $\T^1_X$.
Note that there exists a natural isomorphism
$\T^1_X\stackrel{\sim}{\rightarrow}\Ext^1_{\O_X}(\Omega^1_{X/k},\O_X)$.

Suppose $X$ has a global NCD embedding $X\hookrightarrow V$.
Then by Lemma \ref{tangentloc}, the restriction of the normal bundle
$\N_{X|V}$ to the singular locus $D$ is isomorphic to $\T^1_X$.
Hence we have the following:

\begin{pro}\label{suffcondemb}
If a normal crossing variety $X$ over $k$ is embedded into a smooth $k$-variety
as a normal crossing divisor, then there exists a line bundle $\L_X$ on $X$
such that $\L_X\otimes_{\O_X}\O_D\stackrel{\sim}{\rightarrow}\T^1_X$.
\end{pro}

Let $\Xg\rightarrow\Delta$ be a semistable reduction of schemes,
{\it i.e.},
a flat and generically smooth morphism between regular schemes with
$\Delta$ one-dimensional and every closed fiber is a normal crossing
variety.
Suppose $X\rightarrow \Spec k$ is isomorphic to a closed fiber of this
family.
Then one sees that the normal bundle $\N_{X|V}$ is trivial on $X$, and so
is $\T^1_X$.

\begin{dfn}{\rm (cf. \cite{Fri1})}\ \label{dsemistable}{\rm
A normal crossing variety $X$ is said to be {\it $d$-semistable}
if $\T^1_X$ is the trivial line bundle on $D$.}
\end{dfn}

Due to the above observation, we have the following:

\begin{pro}{\rm (cf. \cite{Fri1})}\label{suffcondsmoothing}
The $d$-semistablilty is a necessary condition for the existence of global
smoothings of $X$.
\end{pro}
\section{Logarithmic embeddings}
In this section, we define the logarithmic embedding of a
normal crossing varieties
(cf. \cite{Ste1}).
This concept is defined in terms of log geometry of Fontaine, Illusie,
and Kazuya Kato (cf. \cite{Kat1}).

Let $X$ be a normal crossing variety over a field $k$.
Suppose that $X$ has a NCD embedding $\iota\colon X\hookrightarrow V$.
We denote the open immersion $V\setminus X\hookrightarrow V$ by $j$.
We define a log structure on $X$ by
$$
\iota^*(\O_V\bigcap j_*\Oi_{V\setminus X})\longrightarrow\O_X,
$$
where $\iota^*$ denotes the pull--back of log structures
(cf. \cite[(1.4)]{Kat1}).
We call this the log structure associated to the NCD embedding
$\iota\colon X\hookrightarrow V$.
For a general normal crossing variety $X$, we cannot define the log structure
of this type on $X$, because $X$ may not have a NCD embedding.
But, as we have seen in Remark \ref{ncvemb}, $X$ has \'{e}tale locally
a NCD embedding.
Then we can consider the log strcuture of this type for a general $X$
defined as follows:

\begin{dfn}{\rm (cf. \cite{Ste1})}\ \label{logembdef}{\rm
A log structure $\M_X\rightarrow\O_X$ is said to be of {\it embedding type},
if the following condition is satisfied:
There exists an \'{e}tale covering $\{\varphi_{\la}\colon U_{\la}\rightarrow
X\}_{\la\in\La}$ by local charts --- with the NCD embeddings
$\iota_{\la}\colon U_{\la}\hookrightarrow V_{\la}$ as in
Definition \ref{ncvchart} --- such that, for each $\la\in\La$, the
restriction
$$
\M_{U_{\la}}\colon =\varphi^*_{\la}\M_X\longrightarrow\O_{U_{\la}}
$$
is isomorphic to the log structure
associated to the NCD embedding $\iota_{\la}$.
If $\M_X\rightarrow\O_X$ is a log structure of embedding type of $X$, we call
the log scheme $(X,\M_X)$ the {\it logarithmic embedding}.
}
\end{dfn}

Let $(X,\M_X)$ be a logarithmic embedding.
We can explicitly write this log structure $\M_X$ \'{e}tale locally.
Let $\nu\colon \widetilde{X}\rightarrow X$ be a normalization of $X$.
Take a local chart $\varphi\colon U\rightarrow X$ with parameters
$z_1,\ldots,z_l$
such that $\M_U\colon =\varphi^*\M_X\rightarrow\O_U$ is the log structure
associated to the NCD embedding $\iota\colon U\hookrightarrow V$.
Let $U=\bigcup^l_{i=1}U_i$ be the decomposition into irreduclbile components,
where $U_i$ is the irreducible component corresponding to the ideal $(z_i)$.
The normalization $\widetilde{U}=\coprod^{l}_{i=1}U_i\rightarrow U$ is denoted
by $\nu_U$.
Note that, due to Lemma \ref{normalization}, we have $U\times_X\widetilde{X}
\cong\widetilde{U}$.
Define a homomorphism of monoids
\begin{equation}\label{logembchart}
\alpha\colon (\nu_U)_*\Na_{\widetilde{U}}\rightarrow\O_U
\end{equation}
by
$\alpha(e_{U_i})=z_i$ for $i=1,\ldots,l$, where $(e_{U_i})$ is the standard
base of $(\nu_U)_*\Na_{\widetilde{U}}=\bigoplus^l_{i=1}\Na_{U_i}$.
Then $\alpha$ induces a log structure
\begin{equation}\label{logembloc}
\Oi_U\bigoplus(\nu_U)_*\Na_{\widetilde{U}}\rightarrow\O_U.
\end{equation}
\begin{pro}\label{logembloc-imp}
The log structure $\M_U\rightarrow\O_U$ is isomorphic to (\ref{logembloc}).
\end{pro}

\pf
Let $Z_1,\ldots,Z_l\in\Gamma(V,\O_V)$ be as in Definition \ref{ncvchart}.
By definition of the log structure associated to the embedding
$\iota\colon U\hookrightarrow V$, these $Z_1,\ldots,Z_l$ are sections of the
sheaf
$\M_U$.
Define a morphism
$$
\psi\colon (\nu_U)_*\Na_{\widetilde{U}}\longrightarrow\M_U
$$
by $\psi(e_{U_i})\colon =Z_i$ for $1\leq i\leq l$.
Let
$$
\widetilde{\psi}\colon (\nu_U)_*\Na_{\widetilde{U}}\rightarrow
\M_U/\Oi_U
$$
be the composition of $\psi$ followed by the natural projection
$\M_U\rightarrow\M_U/\Oi_U$.
It is easy to see that $\widetilde{\psi}$ is injective.
Since sections of $\O_V\cap j_*\Oi_{V\setminus U}$ are precisely those of
$\O_V$
which may take zeros along $\iota(U)$, these are written in the
form $uZ^{a_1}_1\cdots Z^{a_l}_l$ where $u\in\Oi_V$ and $a_1,\ldots,a_l
\in\Na$.
This implies that the morphism $\widetilde{\psi}$ is an isomorphism.
Then, consider the exact sequence of sheaves of monoids
$$
1\rightarrow\Oi_U\rightarrow\M_U\rightarrow\M_U/\Oi_U\rightarrow 1,
$$
where the second arrow is injective.
This exact sequence splits since $\widetilde{\psi}$ is an isomorphism and
$\psi$ defines a cross section $\M_U/\Oi_U\rightarrow\M_U$.
By this, we can easily obtain the desired result.
\qed

Thus, a log structure of embedding type is determined by the morphism
$\alpha\colon (\nu_U)_*\Na_{\widetilde{U}}\rightarrow\O_U$ such that
$\alpha(e_{U_i})$ is a local defining function of the component $U_i$ for
each $i=1,\ldots,l$.
Let ${\alpha}'$ be another such homomorphism.
Then --- replacing $U$ by sufficiently small Zariski open subset ---
we can take $u_i\in\Gamma(U,\Oi_X)$ such that
${\alpha}'(e_{U_i})=u_i\alpha(e_{U_i})$ for each $i$.
Then the isomorphism of log structures of embedding type determined by
$\alpha$ and ${\alpha}'$ is described by the following commutative diagram
\begin{equation}\label{logemb-iso}
\begin{array}{ccccc}
\Oi_U\oplus(\nu_U)_*\Na_{\widetilde{U}}&&\stackrel{\phi}{\longrightarrow}&&
\Oi_U\oplus(\nu_U)_*\Na_{\widetilde{U}}\\
&\llap{$\mbox{\rm by}\ {\alpha}'$}\searrow&&\swarrow\rlap{$\mbox{\rm by}\
\alpha$}\\
&&\O_X\rlap{,}
\end{array}
\end{equation}
where $\phi$ is defined by $\phi(1,e_{U_i})=(u_i,e_{U_i})$ for each
$i=1,\ldots,l$.
In particular, the log structure of embedding type exists \'{e}tale locally,
and is unique up to isomorphisms.

\begin{cor}
For any logarithmic embedding $(X,\M_X)$,
we have an exact sequence of abelian sheaves
\begin{equation}\label{abelian}
1\longrightarrow\Oi_X\longrightarrow\gp{\M_X}\longrightarrow
\nu_*\Z_{\widetilde{X}}\longrightarrow 0.
\end{equation}
\end{cor}

\pf
Due to the local expression (\ref{logembloc}).
\qed

In the rest of this section, we prove the following theorem, which is the
main theorem of this paper.

\begin{thm}\label{mainthm}
For a normal crossing variety $X$, the logarithmic embedding of $X$ exists
if and only if there exists a line bundle ${\cal L}$ on $X$ such that
${\cal L}\otimes_{\O_X}\O_D\stackrel{\sim}{\rightarrow} \T^1_X$.
\end{thm}

For the proof of this theorem, we shall prove some lemmas as follows.
Let $(\varphi\colon U=\Spec A\rightarrow X; z_1,\ldots,z_l)$ be a local chart
on $X$.
Let the NCD embedding $U\hookrightarrow V=\Spec R$ and the ideals $I_j,I,
J_j,J$ of $R$ be as in the previous section.

\begin{lem}\label{lem-1}
The natural morphism
$$
\bigoplus^l_{j=1}J_j/I\longrightarrow J/I
$$
of $A$-modules, induced by $J_j\hookrightarrow J$, is an isomorphism.
\end{lem}

\pf
The surjectivity is clear.
We are going to show the injectivity.
Take $a_jZ_1\cdots\widehat{Z_j}\cdots Z_l\in J_j$ --- where $Z_1,\ldots,Z_l$
are as in the previous section --- for $1\leq j\leq l$ such that
$$
\sum^l_{j=1}a_jZ_1\cdots\widehat{Z_j}\cdots Z_l=b\cdot Z_1\cdots Z_l,
$$
where $a_j,b\in R$.
Since $R$ is an integral domain, $a_j$ is divisible by $Z_j$, and hence, we
have $a_jZ_1\cdots\widehat{Z_j}\cdots Z_l\equiv 0\ (\mbox{\rm mod}\,I)$.
\qed

Let $\pi_j\colon I_j/II_j\rightarrow I_j/I$ and $q_j\colon I_j/I
\rightarrow I_j/JI_j(\cong I_j/II_j\otimes_AQ \ \mbox{\rm where}\ Q=R/J)$ be
the
natural projections and set $p_j\colon=q_j\circ\pi_j$.
Let $q\colon I/I^2\rightarrow I/JI(\cong I/I^2\otimes_AQ)$ be the natural
projection.

\begin{lem}\label{lem-2}
Let $M_1,\ldots,M_l$ be free $A$-modules of rank one and set
$M\colon=M_1\otimes_A\cdots\otimes_AM_l$.
Suppose we are given an $A$-module isomorphism
$\widetilde{g}\colon M\stackrel{\sim}{\rightarrow}I/I^2$ and $A$-module
homomorphisms
$g_j\colon M_j\rightarrow I_j/I$, for $1\leq j\leq l$, such that,
\begin{enumerate}
\item
for each $j$, there exists a free generator $\delta_j$ of $M_j$ such that
$g_j(\delta_j)=z_j$,
\item
$(q_1\circ g_1)\otimes_Q\cdots\otimes_Q(q_l\circ g_l)=q\circ\widetilde{g}$.
\end{enumerate}
Then there exists a unique collection $\{\widetilde{g}_j\colon M_j
\stackrel{\sim}{\rightarrow}I_j/II_j\}^{l}_{j=1}$ of $A$-isomorphisms
such that
$\pi_j\circ\widetilde{g}_j=g_j$ for each $j$ and
$\widetilde{g}_1\otimes_A\cdots\otimes_A\widetilde{g}_l=\widetilde{g}$.
\end{lem}

\pf
We fix the free generators $\delta_j$ of $M_j$ as above.
Then $M$ is generated by $\delta_1\otimes\cdots\otimes\delta_l$.
Set $\widetilde{g}(\delta_1\otimes\cdots\otimes\delta_l)=
v\zeta_1\otimes\cdots\otimes\zeta_l$ where $v\in A^{\times}$.
By the second condition, we have
$v\equiv 1\ (\mbox{\rm mod}\,J/I)$, {\it i.e.},
$$
v=1+\sum^l_{j=1}a_jz_1\cdots\widehat{z_j}\cdots z_l
$$
for $a_j\in A$.
We set $u_j=1+a_jz_1\cdots\widehat{z_j}\cdots z_l$ and define
$\widetilde{g}_j$ by $\widetilde{g}_j(\delta_j)\colon=u_j\zeta_j$ for
$1\leq j\leq l$.
Then, since $v=u_1\cdots u_l$, each $u_j$ is a unit in $A$ and
$\widetilde{g}_j$ is an isomorphism. Moreover, we have
$\widetilde{g}_1\otimes\cdots\otimes\widetilde{g}_l=\widetilde{g}$
as desired.
The uniqueness follows from Lemma \ref{lem-1}.
\qed

\vspace{3mm}
{\sc Proof of Theorem \ref{mainthm}}.
We first prove the ``if'' part. This part is divided into four steps.

{\sc Step 1}:
Here, we shall describe the log structure of embedding type by another
\'{e}tale local expression.
Let $(\varphi\colon U=\Spec A\rightarrow X; z_1,\ldots,z_l)$ be a local chart.
For $m=(m_1,\ldots,m_l)\in\Na^l$, define an $A$-module $P_m$ by
$$
P_m\colon=(I_1/II_1)^{\otimes m_1}\otimes_A\cdots\otimes_A
(I_l/II_l)^{\otimes m_l}.
$$
Each $P_m$ is a free $A$-module of rank one and $P_{(1,\ldots,1)}\cong
I/I^2$.
The natural projections $\pi_j$ induce a natural $A$-homomorphism
$$
\sigma_m\colon P_m\longrightarrow A.
$$
Define a monoid
$$
M\colon =\left\{
\begin{array}{c|l}
(m,a)&m\in\Na^l,\\
&a:\mbox{a generator of $P_m$}
\end{array}
\right\},
$$
and a homomorphism $M\rightarrow A$ of monoids by $(m,a)\mapsto\sigma_m(a)$.
Then the associated log structure $\alpha_U\colon\M_U\rightarrow\O_U$ of the
pre--log structure $M\rightarrow A$ is that of embedding type on $U$.

{\sc Step 2}:
Now, we assume that we are given a line bundle $\L$ on $X$ satisfying
$\L\otimes_{\O_X}\O_D\cong(\T^1_X)^{\vee}$.
Suppose we have two local charts $(\varphi\colon U\rightarrow X;
z_1,\ldots,z_{l})$ and $(\varphi'\colon U'\rightarrow X;
z'_1,\ldots,z'_{l'})$ and an \'{e}tale morphism $\psi\colon U\rightarrow U'$
such that $\varphi=\varphi'\circ\psi$.
For these local charts, we use the notation as in the previous
section; such as $U=\Spec A\hookrightarrow V=\Spec R$,
$U'=\Spec A'\hookrightarrow V'=\Spec R'$, $f\colon A'\rightarrow A$,
$I$, $I'$, etc.
As in the previous section, we may assume
$(I'_j/I')\otimes_{A'}A=I_j/I$ as ideals in $A$ for
$1\leq j\leq l$, and set $f(z'_j)=u_jz_j$ (each $u_j$ is determined up to
modulo $J_j/I$).
To give the line bundle $\L$ as above is equivalent to give a compatible
system of isomorphisms
$$
\widetilde{\tau}\colon
I'/I'^2\otimes_{A'}A\stackrel{\sim}{\longrightarrow}I/I^2,
$$
for all such $U\rightarrow U'$, with $\widetilde{\tau}\otimes_AQ=\tau$,
where $\tau$ is defined as in (\ref{def-tau2}).
Then we shall show that $\widetilde{\tau}$ induces canonically an isomorphism
of log structures $\psi^*\M_{U'}\stackrel{\sim}{\rightarrow}\M_U$, and
prove that these
isomorphisms form so a compatible system that the log structures $\M_U$ glue to
a log structure of embedding type on $X$.
Moreover --- since local charts form an \'{e}tale open basis (Corollary
\ref{specialcov}) --- we can pass through this procedure replacing $U$ by
its Zariski open subset if necessary.
In particular, we may assume that each $u_j$ as above is a unit in $A$,
because $(u_j\,\mbox{\rm mod}\,J/I)$ is a unit in $A/(J/I)$ (in case $l>1$).
Fix a locally constant section $w\in\H^0(D,\Oi_D)$.
(Actually, we can take $w$ as any
global section in $\H^0(D,\Oi_D)$ but, if we do so, the following argument
have to be modified slightly.)

{\sc Step 3}:
(i) If $l=l'=1$, {\it i.e.}, $I_1=I$ and $I'_1=I'$, then we set
$\widetilde{\tau}_1\colon I'_1/I'I'_1\otimes_{A'}A\stackrel{\sim}{\rightarrow}
I_1/II_1$ by $\widetilde{\tau}_1\colon=\widetilde{\tau}$.

(ii) If $l=1$ and $l'>1$, we define
$\widetilde{\tau}_1\colon I'_1/I'I'_1\otimes_{A'}A\stackrel{\sim}{\rightarrow}
I_1/II_1$ as follows:
Suppose $\widetilde{\tau}$ maps $\zeta'_1\otimes\cdots\otimes\zeta'_{l'}\otimes
1$ to
$v\zeta_1$, where $v\in A^{\times}$.
Let $\widetilde{\rho}_j\colon I'_i/I'I'_i\otimes_{A'}A\rightarrow A$ be as
(\ref{def-proj}), for $1\leq i\leq l'$.
Suppose, moreover, each $\widetilde{\rho}_i$, for $i>1$, maps $\zeta'_i\otimes
1$ to $v_i\in A^{\times}$.
Then, define $\widetilde{\tau}_1$ by $\widetilde{\tau}_1(\zeta'_1\otimes 1)
\colon=w_Uvv^{-1}_2\cdots v^{-1}_{l'}\zeta_1$, where $w_U$ is a non--zero
scalar which coincides with $w$ restricted to $D_U$.

(iii) Suppose $l>1$ and $l'>1$.
We claim that, under the conditions
\begin{equation}\label{cond1}
\pi_j\circ\widetilde{\tau}_j=\widetilde{\rho}_j,\ (1\leq j\leq l)
\end{equation}
and
\begin{equation}\label{cond2}
\widetilde{\tau}_1\otimes_A\cdots\otimes_A\widetilde{\tau}_l\otimes_A
\widetilde{\rho}_{l+1}\otimes_A\cdots\otimes_A\widetilde{\rho}_{l'}=
\widetilde{\tau},
\end{equation}
the $A$-isomorphisms
$$
\widetilde{\tau}_j\colon I'_j/I'I'_j\otimes_{A'}A\stackrel{\sim}
{\longrightarrow}I_j/II_j
$$
exist uniquely for $1\leq j\leq l$.
Set $M_j\colon=I'_j/I'I'_j\otimes_{A'}A$ and
$g_j\colon=\widetilde{\rho}_j$ for $1\leq j\leq l$.
Define
$\widetilde{g}$ by $\widetilde{g}\otimes_A\widetilde{\rho}_{l+1}\otimes_A
\cdots\otimes_A\widetilde{\rho}_{l'}=\widetilde{\tau}$ (this is possible
since $\widetilde{\rho}_i(\zeta'_i\otimes 1)$ is a unit element in $A$ for
$i>l$), which is obviously an isomorphism.
Then --- since we assumed each $u_j$ to be a unit in $A$ --- $M_j\colon=
I'_j/I'_jI'\otimes_{A'}A$, $g$, and $g_j$ satisfy the conditions in Lemma
\ref{lem-2}.
Hence our claim follows from this lemma.

Note that, in any cases, we have the following commutative diagram:
\begin{equation}\label{com-mor}
\begin{array}{ccc}
I'_j/I'I'_j&\longrightarrow&I_j/II_j\\
\llap{$\pi'_j$}\Bigdownarrow&&\Bigdownarrow\rlap{$\pi_j$}\\
A'&\underrel{\longrightarrow}{f}&A\rlap{,}
\end{array}
\end{equation}
for $1\leq j\leq l$; this follows from (\ref{cond1}) in case $l, l'>1$, and is
quite obvious in the other cases.

{\sc Step 4}:
These morphisms $\widetilde{\tau}_j$ induce the morphisms
$$
\gamma_{m'}\colon P'_{m'}\longrightarrow P_{m},
$$
where $m=(m_1,\ldots,m_l)$ for $m'=(m_1,\ldots,m_{l'})\in\Na^{l'}$.
Then these $\gamma_{m'}$ induce naturally a morphism of monoids
$M'\rightarrow M$ compatible with $M'\rightarrow A'$, $M\rightarrow A$ and
$f$.
By the construction of these morphisms, the induced morphism of sheaves
of monoids
$\gamma\colon \psi^*\M_{U'}\stackrel{\sim}{\rightarrow}\M_U$ is an isomorphism.
By the commutative diagram (\ref{com-mor}), this isomorphism commutes the
following diagram:
$$
\begin{array}{ccc}
\psi^*\M_{U'}&\stackrel{\sim}{\longrightarrow}&\M_{U}\\
\llap{$\psi^*\alpha_{U'}$}\Bigdownarrow&&\Bigdownarrow\rlap{$\alpha_U$}\\
\O_U&=&\O_U\rlap{;}
\end{array}
$$
hence $\gamma$ is an isomorphism of log structures.
Our construction of the isomorphism $\gamma$ is canonical in the following
sense: Suppose we are given a sequence of \'{e}tale morphisms
$U\stackrel{\psi}{\rightarrow}U'\stackrel{\psi'}{\rightarrow}U''$ of
local charts (with $U$ and $U'$ sufficiently small),
we have $\gamma''=\gamma\circ(\psi^*\gamma')$, where
$\gamma\colon \psi^*\M_{U'}\stackrel{\sim}{\rightarrow}\M_U$,
$\gamma'\colon \psi'^*\M_{U''}\stackrel{\sim}{\rightarrow}\M_{U'}$ and
$\gamma''\colon \psi^*\psi'^*\M_{U''}\stackrel{\sim}{\rightarrow}\M_U$
are the isomorphisms of log structures defined as above corresponding to
$\psi$, $\psi'$ and $\psi'\circ\psi$, respectively.
This follows from the naturality of $\pi_j$ and $\widetilde{\rho}_j$, and
the compatibility of $\widetilde{\tau}$'s. Then there exists a unique
log structure $\M_X$ on $X$ which is of embedding type.
Hence the ``if'' part is now proved.

Conversely, suppose we are given a log structure $\M_X$ of embedding type.
Then we have an exact sequence (\ref{abelian}) of abelian sheaves.
Considering the cohomology exact sequence, we obtain a morphism
$$
\delta\colon \H^0(X,\nu_*\Z_{\widetilde{X}})\longrightarrow\H^1(X,\Oi_X)
(\cong\Pic X).
$$
In $\H^0(X,\nu_*\Z_{\widetilde{X}})$, we consider the element $\diag$
which is defined by
the image of $1\in\Z_X$ under the diagonal morphism $\Z_X\rightarrow
\nu_*\Z_{\widetilde{X}}$.
Then $\delta(\diag)$ defines a line bundle $\L=\L_{\M_X}$ on $X$.
We shall show that this line bundle satisfies $\L\otimes_{\O_X}\O_D
\stackrel{\sim}{\rightarrow}(\T^1_X)^{\vee}$.

The line bundle $\L$ is constructed as follows: the inverse image of
$\diag$ under $\gp{\M_X}\rightarrow\nu_*\Z_{\widetilde{X}}$ defines a
principally homogeneous space over $\Oi_X$ and hence defines a line bundle,
which is nothing but $\L$.
Let $U=\Spec A$ be a local chart as above.
Then the inverse image of $\diag$ restricted to $U$ gives a generator
of an $A$-module $I/I^2$ which is --- due to Lemma \ref{tangentloc} ---
a local lifting of $\T^1_X$ restricted to
$U$. Hence $\L$ satisfies the desired condition.
\qed

\begin{rem}\label{important}{\rm
1.
As we have seen above, the log structure of embedding type exists locally and
is unique up to isomorphisms.
The sheaf of germs of automorphisms of such a log structure
is naturally isomorphic to $\K$, where $\K$ is defined by the exact sequence
\begin{equation}\label{important2}
1\longrightarrow\K\longrightarrow\Oi_X\longrightarrow
\Oi_D\longrightarrow 1.
\end{equation}
This can be shown by the following steps:
(i) any automorphism over a sufficiently small local chart $U$ is given by $
\phi$ in the diagram
(\ref{logemb-iso}) with $\alpha=\alpha'$; (ii) $\phi$ is determined by
$\{u_i\}$ with $u_i\in\Gamma(U,\Oi_X)$ such that $z_i=u_i\cdot z_i$ for each
$i$; (iii) hence such $u_i$'s are written in the form of $u_i=
1+a_i\cdot z_1\cdots\widehat{z_i}\cdots z_l$; (iv) due to Lemma \ref{lem-1},
to give a system
$\{u_i\}$ is equivarent to give $u=u_1\cdots u_l$ which is a section of $\K$.
Hence the obstruction for the existence of log structures of embedding type
lies
in $\H^2(X,\K)$.
The proof of Theorem \ref{mainthm} shows that this class coincides with the
obstruction class for a lifting of $(T^1_X)^{\vee}$ on $X$, {\it i.e.}, the
image
of $(T^1_X)^{\vee}$ under $\H^1(D,\Oi_D)\rightarrow\H^2(X,\K)$.

2.
One sees easily --- by the proof of Theorem \ref{mainthm} --- that there
exists a natural surjective map
\begin{equation}\label{important1}
\left\{
\begin{array}{ll}
\mbox{isom. class of log structures}\\
\mbox{of embedding type on $X$}
\end{array}
\right\}
{\longrightarrow}
\left\{
\begin{array}{c|c}
\L\in\Pic X&\L\otimes_{\O_X}\O_D\stackrel{\sim}{\rightarrow}(\T^1_X)^{\vee}
\end{array}
\right\}
\end{equation}
by $\M\mapsto\L_{\M}$, where $\L_{\M}$ is defined as in the proof of Theorem
\ref{mainthm}.
If $\M_X$ is associated to a global NCD embedding $X\hookrightarrow V$, then
$\L_{\M_X}$ is nothing but the conormal bundle of $X$ in $V$.
The set of isomorphism classes of log structures of embedding type on $X$, is
a principally homogeneous space over $\H^1(X,\K)$.
Then one sees easily that the map (\ref{important1}) is equivariant to
$\H^1(X,\K)\rightarrow
\Ker(\H^1(X,\Oi_X)\rightarrow\H^1(D,\Oi_D))$ induced by the
cohomology exact sequence of (\ref{important2}). In particular, if $X$ is
proper and $D$ is connected, the map (\ref{important1}) is a bijection
since $\H^1(X,\K)\stackrel{\sim}{\rightarrow}\Ker(\H^1(X,\Oi_X)\rightarrow
\H^1(D,\Oi_D))$; in this case, the logarithmic embeddings are determined by
their ``normal bundles.''

3.
By the exact sequence (\ref{abelian}), a log structure of embedding type
$\M_X$ on
$X$ defines an extension class in
$\Extr^1_{\Z_X}(\nu_*\Z_{\widetilde{X}},\Oi_X)$.
Under the morphism $\Extr^1_{\Z_X}(\nu_*\Z_{\widetilde{X}},\Oi_X)
\rightarrow\Extr^1_{\Z_X}(\Z_X,\Oi_X)$, induced by the diagonal morphism
$\Z_X\rightarrow\nu_*\Z_{\widetilde{X}}$, and the natural identification
$\Extr^1_{\Z_X}(\Z_X,\Oi_X)\stackrel{\sim}{\rightarrow}\Pic X$, this class is
mapped
to the class corresponding to the line bundle $\L_{\M_X}$ defined as above.
(The proof is straightforward and left to the reader.)
}
\end{rem}
\section{Logarithmic semistable reductions}

\begin{dfn}{\rm (cf. \cite{Kaj1}, \cite{K-N1})}\label{logsemidef}\ {\rm
A log strcuture of embedding type $\M_X\rightarrow\O_X$ is said to be
of {\it semistable type}, if there exists a homomorphism
$\Z_X\rightarrow\gp{\M_X}$ of abelian sheaves on $X$ such that the diagram
$$
\begin{array}{ccccc}
\gp{\M_X}&&\longrightarrow&&\nu_*\Z_{\widetilde{X}}\\
&\nwarrow&&\nearrow\rlap{$\diag$}\\
&&\Z_X
\end{array}
$$
commutes, where $\diag\colon\Z_X{\rightarrow}\nu_*\Z_{\widetilde{X}}$ is
the diagonal homomorphism, and $\gp{\M_X}\rightarrow\nu_*\Z_{\widetilde{X}}$
is the projection in (\ref{abelian}).
}
\end{dfn}

If $\M_X$ is a log strcuture of semistable type, the homomorphism
$\Z_X\rightarrow\gp{\M_X}$ induces the homomorphism $\Na_X\rightarrow\M_X$
of monoids by the following Cartesian diagram:
$$
\begin{array}{ccc}
\Na_X&\longrightarrow&\M_X\\
\Bigdownarrow&&\Bigdownarrow\\
\Z_X&\longrightarrow&\gp{\M_X}\rlap{;}
\end{array}
$$
this follows easily from the local expression (\ref{logembloc}).
This morphism defines a morphism of log schemes
$$
(X,\M_X)\longrightarrow(\Spec k,\Na)
$$
Here, $(\Spec k,\Na)$ is the {\it standard point} defined by $\Na\rightarrow k$
which maps $m\in\Na$ to $0^m$.
We call this morphism of log schemes the {\it logarithmic semistable
reduction}.
Logarithmic semistable reductions are {\it log smooth}
in the sense of \cite{Kat1}.

\begin{rem}\label{genuine-semistable}{\rm
Let $f\colon\Xg\rightarrow\Delta$ be a semistable reduction of schemes;
{\it i.e.}, a proper flat generically smooth morphism $f$ with $\Xg$ a
regular scheme and $\Delta$ a one-dimensional regular local scheme,
with the closed fiber $X\rightarrow 0=\Spec k$ a normal crossing variety.
Then this morphism induces canonically a logarithmic semistable reduction
$(X,\M_X)\rightarrow(\Spec k,\Na)$ on the closed fiber as follows:
We define a log structure $\M_{\Xg}\rightarrow\O_{\Xg}$ by
$$
\M_{\Xg}\colon=\O_{\Xg}\bigcap j_*\Oi_{\Xg\setminus X}
\longhookrightarrow\O_{\Xg}
$$
where $j\colon \Xg\setminus X\hookrightarrow\Xg$ is an open
immersion.
Take a local parameter $t\in\O_{\Delta}$ around $0=\Spec k$.
Then $f^{-1}(t)$ belongs to $\M_{\Xg}$.
We define a homomorphism of monoids $\Na_{\Xg}\rightarrow\M_{\Xg}$ by
$1\mapsto f^{-1}(t)$. Then this homomorphism extends to a morphism of
log schemes
\begin{equation}\label{france}
(\Xg,\M_{\Xg})\longrightarrow(\Delta,0),
\end{equation}
where the log structure on $\Delta$ is the associated log structure of
$$
\Na\longrightarrow\O_{\Delta}\quad \mbox{\rm by}\quad
m\mapsto t^m.
$$
Taking the pull--back of (\ref{france}) to the closed fiber, we
get a logarithmic semistable reduction.
Note that the monoid morphism $\Na_{\Xg}\rightarrow\M_{\Xg}$ induces
$\Z_X\rightarrow\gp{\M_X}$ which satisfies the condition in Definition
\ref{logsemidef}.
Hence, such a morphism $\Z_X\rightarrow\gp{\M_X}$ for a general log structure
of semistable type can be regarded as a ``parametrization.''
}
\end{rem}

\begin{rem}{\rm
The logarithmic semistable reduction induced by a semistable reduction
family, as in Remark \ref{genuine-semistable}, is regarded as the
``closed fiber'' of the morphism (\ref{france}) of log schemes.
Then, conversely, one can consider the theory of deformations which
deal with liftings of the logarithmic semistable reductions.
This is nothing but the {\it logarithmic deformation} of Kawamata--Namikawa
\cite{K-N1}, and also a part of the {\it log smooth deformation} developed
in \cite{Kat2}.
}
\end{rem}

Using Theorem \ref{mainthm} --- which is proved in the previous section ---
we get a new proof of the theorem of Kawamata--Namikawa as follows:
\begin{thm}\label{mainthm2}{\rm (cf. \cite{K-N1})}\
For a normal crossing variety $X$, the log structure of semistable type
on $X$ exists if and only if $X$ is $d$-semistable.
\end{thm}

To prove the theorem, we need the following lemma:
\begin{lem}\label{l1}
Let $(X,\M_X)$ be a logarithmic embedding.
Consider the exact sequence (\ref{abelian}) of abelian sheaves and the
induced morphism
$$
\Hom_{\Z_X}(\Z_X,\nu_*\Z_{\widetilde{X}})\stackrel{\delta}{\longrightarrow}
\Extr^1_{\Z_X}(\Z_X,\Oi_X).
$$
Let $\diag\in\Hom_{\Z_X}(\Z_X,\nu_*\Z_{\widetilde{X}})$ be the diagonal
morphism.
Then, under the natural identification $\Extr^1_{\Z_X}(\Z_X,\Oi_X)
\stackrel{\sim}{\rightarrow}\Pic X$, we have
$$
\delta(\diag)=[\L_{\M_X}],
$$
where $\L_{\M_X}$ is the line bundle defined in the previous section.
\end{lem}

\pf
This lemma follows from the commutative diagram
$$
\begin{array}{ccc}
\Hom_{\Z_X}(\Z_X,\nu_*\Z_{\widetilde{X}})&\stackrel{\delta}{\longrightarrow}&
\Extr^1_{\Z_X}(\Z_X,\Oi_X)\\
\llap{$\cong$}\Bigdownarrow&&\Bigdownarrow\rlap{$\cong$}\\
\H^0(X,\nu_*\Z_{\widetilde{X}})&\longrightarrow&\H^1(X,\Oi_X)
\end{array}
$$
where the vertical morphisms are natural isomorphisms and the definition
of the line bundle $\L_{\M_X}$.
\qed

\vspace{3mm}
{\sc Proof of Theorem \ref{mainthm2}.}\hspace{2mm}
Suppose $\M_X$ is a log structure of semistable type.
Consider the exact sequence
\begin{equation}\label{extension}
\Hom_{\Z_X}(\Z_X,\gp{\M_X})\stackrel{\pi}{\longrightarrow}
\Hom_{\Z_X}(\Z_X,\nu_*\Z_{\widetilde{X}})\stackrel{\delta}{\longrightarrow}
\Extr^1_{\Z_X}(\Z_X,\Oi_X)
\end{equation}
induced by (\ref{abelian}).
The ``parametrization'' morphism $\Z_X\rightarrow\gp{\M_X}$ is mapped to
$\diag$ by $\pi$.
This implies that the line bundle $\L_{\M_X}$ is trivial.
Then so is $(T^1_X)^{\vee}$ because $\L_{\M_X}\otimes_{\O_X}\O_D$ is
isomorphic to $(T^1_X)^{\vee}$.

Conversely, if $X$ is $d$-semistable, there exists at least one
log structure of embedding type on $X$ due to Theorem \ref{mainthm}.
Since $\L_{\M_X}\otimes_{\O_X}\O_D$ is trivial, we can take the log structure
$\M_X$ of embedding type such that the corresponding line bundle $\L_{\M_X}$
is trivial (due to the natural surjection (\ref{important1})).
Since the obstruction for the existence of a morphism
$\Z_X\rightarrow\gp{\M_X}$ which is mapped to $\diag$, is nothing but the
class $[\L_{\M_X}]$, we deduce that $\M_X$ is of semistable type.
\qed

As is shown in the above proof, the log structure of semistable type on $X$
is --- considering the natural surjection (\ref{important1}) --- the log
structure of embedding type which is mapped
to the trivial bundle on $X$. Hence we have the following:

\begin{cor}
Let $X$ be a proper, $d$-semistable normal crossing variety, and assume that
the singular locus $D$ is connected.
Then, the log structure of semistable type on $X$ exists uniquely.
\end{cor}

\begin{exa}{\rm
Let $X\colon =X_0\cup\cdots\cup X_N$ be a chain of surfaces defined as follows:
Each $X_i$ is the Hirzebruch surface of degree $a_i\leq 0$.
The surfaces $X_{i-1}$
and $X_i$ are connected by identifying the section $s'_{i-1}$ on $X_{i-1}$
and the one $s_i$ on $X_i$, where $(s'_{i-1})^2=a_{i-1}$ and $(s_i)^2=
-a_i$, for $1\leq i\leq N$. Then $X$ has a
log structure of embedding type if and only if $a_i|(a_{i-1}+a_{i+1})$
for $1\leq i\leq N-1$, while $X$ has a log structure of semistable type
if and only if $a_0=a_1=\cdots =a_N$.
}
\end{exa}
%
%
\begin{small}

\end{small}
\end{document}